\documentclass[a4paper,12pt]{article}
\usepackage{amsthm,amsmath}
\usepackage{amssymb}
\usepackage[all]{xypic}
\usepackage{enumerate}
\usepackage{a4wide}
\usepackage{hyperref}

\swapnumbers
\newtheorem{theorem}{Theorem}[section]
\newtheorem{proposition}[theorem]{Proposition}
\newtheorem{corollary}[theorem]{Corollary}
\newtheorem{lemma}[theorem]{Lemma}
\newtheorem*{theorem*}{Theorem}
\newtheorem*{proposition*}{Proposition}
\newtheorem*{corollary*}{Corollary}
\newtheorem*{lemma*}{Lemma}
\theoremstyle{definition}
\newtheorem{definition}[theorem]{Definition}
\newtheorem{punto}[theorem]{}
\newtheorem{example}[theorem]{Example}
\newtheorem{remark}[theorem]{Remark}
\newtheorem*{remark*}{Remark}

\newtheorem*{definition*}{Definition}

\newcommand{\coring}[1]{\mathfrak{#1}}
\newcommand{\tensor}[1]{\otimes_{#1}}
\newcommand{\tensfun}[1]{\underset{{#1}}{\otimes}}
\newcommand{\tensdia}[1]{\otimes_{#1}}
\newcommand{\rcomod}[1]{\mathcal{M}^{#1}}
\newcommand{\rmod}[1]{\mathcal{M}_{#1}}
\newcommand{\lmod}[1]{{}_{#1}\mathcal{M}}
\newcommand{\bcomod}[2]{{}^{#1}\mathcal{M}^{#2}}
\newcommand{\cotensor}[1]{\square_{#1}}
\newcommand{\lcomod}[1]{{}^{#1}\mathcal{M}}

\renewcommand{\hom}[3]{\mathrm{Hom}_{#1}(#2,#3)}

\begin{document}
\title{Coseparable corings}
\author{J. G\'omez-Torrecillas \\
\normalsize Departamento de \'{A}lgebra \\
\normalsize Facultad de Ciencias \\
\normalsize Universidad de Granada\\
\normalsize E18071 Granada, Spain \\
\normalsize e-mail: \textsf{torrecil@ugr.es} \and A. Louly \\
\normalsize D\'{e}partment de Math\'{e}matiques \\
\normalsize Facult\'{e} de Sciences \\
\normalsize Universit\'{e} Abdelmalek Esaadi \\
\normalsize B.P. 2121 T\'{e}touan, Morocco }

\date{}

\maketitle

\section*{Introduction}

M. Sweedler \cite{Sweedler:1975} introduced the notion of a coring
as a generalization of the concept of coalgebra in order to study
the intermediate division rings for an extension of division
rings. It turns out that this formalism embodies several kinds of
relative module categories. Thus, graded modules, Doi-Hopf modules
and, more generally, entwined modules are instances of comodules
over suitable corings (see \cite{Brzezinski:2000unp} and its
references). This explains the increasing interest in the study of
corings, as one could expect that the development of the theory of
corings will contribute significantly to the understanding of the
aforementioned relative modules.

Coseparable corings were investigated from a (co)homological point
of view by F. Guzman \cite{Guzman:1989}, dualizing earlier results
by M. Kleiner in the framework of rings \cite{Kleiner:1985}, and
extending Y. Doi's theory on coseparable coalgebras
\cite{Doi:1981}. In this paper, we will relate coseparability and
semisimplicity in the framework of coring theory, proving the
following theorem.

\begin{theorem*}
Let $A$ be a separable algebra over a field $k$. The following
conditions are equivalent for an $A$--coring $\coring{C}$.
\begin{enumerate}[(i)]
\item $\coring{C}$ is coseparable;
\item $\coring{C} \tensor{k} K$ is a semisimple $A \tensor{k} K$--coring for every field
extension $k \subseteq K$;
\item $\coring{C} \tensor{k} \coring{D}$ is a semisimple $A
\tensor{k} B$--coring for every semisimple $B$--coring
$\coring{D}$ ($B$ is any $k$--algebra);
\item $\coring{C} \tensor{k} \coring{C}^{\circ}$ is a semisimple $A \tensor{k} A^{\circ}$--coring.
\end{enumerate}
\end{theorem*}

 When particularized to coalgebras over fields, this theorem gives
 \cite[Proposition 12]{Doi:1981} and \cite[Theorem
 3.4]{Castano/Gomez/Nastasescu:1997}. The notions involved in the statements of our theorem, like
 semisimple coring, opposite coring, or tensor product of corings, are defined
 below in this note.

\section{Bicomodules}\label{bicomodules}

We first recall from \cite{Sweedler:1975} the notion of a coring.
The concepts of comodule and bicomodule over a coring are
generalizations of the corresponding notions for coalgebras, and
were considered in \cite{Guzman:1989}.

Throughout this paper, $A, A', \dots, B, \dots$ denote associative
and unitary algebras over a commutative ring $k$. The tensor
product over $A$ is denoted by $\tensor{A}$. We shall sometimes
replace $\tensor{k}$ by $\tensor{}$.

\begin{punto}\label{coringdef}\textbf{Corings.}
 An $A$--\emph{coring} is a three-tuple
$(\coring{C},\Delta_{\coring{C}},\epsilon_{\coring{C}})$
consisting of an $A$--bimodule $\coring{C}$ and two $A$--bimodule
maps
\begin{equation}
\xymatrix{ \Delta_{\coring{C}} : \coring{C} \ar[r] & \coring{C}
\tensdia{A} \coring{C} & \epsilon_{\coring{C}} : \coring{C} \ar[r]
& A}
\end{equation}
such that the diagrams
\begin{equation}
\xymatrix{ \coring{C} \ar^{\Delta_{\coring{C}}}[rr]
\ar_{\Delta_{\coring{C}}}[d]  & & \coring{C} \tensdia{A}
\coring{C}
\ar^{\coring{C} \tensfun{A} \Delta_{\coring{C}}}[d] \\
\coring{C} \tensdia{A} \coring{C} \ar^{\Delta_{\coring{C}}
\tensfun{A} \coring{C}}[rr] & & \coring{C} \tensdia{A} \coring{C}
\tensdia{A} \coring{C}}
\end{equation}
and
\begin{equation}
\xymatrix{ \coring{C} \ar^{\Delta_{\coring{C}}}[r] \ar_{\cong}[dr]
& \coring{C} \tensdia{A} \coring{C} \ar^{\coring{C} \tensfun{A}
\epsilon_{\coring{C}}}[d] \\
& \coring{C} \tensdia{A} A } \qquad \xymatrix{ \coring{C}
\ar^{\Delta_{\coring{C}}}[r] \ar_{\cong}[dr] & \coring{C}
\tensdia{A} \coring{C}
\ar^{\epsilon_{\coring{C}} \tensfun{A} \coring{C}}[d] \\
& A \tensdia{A} \coring{C} }
\end{equation}
commute.
\end{punto}

From now on, $\coring{C}, \coring{C}', \dots, \coring{D}, \dots$
will denote corings over $A, A', \dots, B,\dots $, respectively.

\begin{punto}\label{comoddef}\textbf{Comodules.}
A \emph{right $\coring{C}$--comodule} is a pair $(M, \rho_M)$
consisting of a right $A$--module $M$ and an $A$--linear map
$\rho_M : M \rightarrow M \tensor{A} \coring{C}$ such that the
diagrams
\begin{equation}
\xymatrix{ M \ar^{\rho_M}[rr] \ar_{\rho_M}[d] & {} & M \tensdia{A}
\coring{C} \ar^{M \tensfun{A} \Delta_{\coring{C}}}[d] \\
M \tensdia{A} \coring{C} \ar^{\rho_M \tensfun{A} \coring{C}}[rr] &
{} & M \tensdia{A} \coring{C} \tensdia{A} \coring{C} } \qquad
\xymatrix{ M \ar^{\rho_M}[r] \ar_{\cong}[dr] & M \tensdia{A}
\coring{C}
\ar^{M \tensfun{A} \epsilon_{\coring{C}}}[d] \\
& M \tensdia{A} A }
\end{equation}
commute. Left $\coring{C}$--comodules are similarly defined; we
use the notation $\lambda_M$ for their structure maps. A
\emph{morphism} of right $\coring{C}$--comodules $(M,\rho_M)$ and
$(N, \rho_N)$ is an $A$--linear map $f : M \rightarrow N$ such
that the following diagram is commutative.
\begin{equation}
\xymatrix{
M \ar^f[rr] \ar^{\rho_M}[d] & & N \ar^{\rho_N}[d] \\
M \tensdia{A} \coring{C} \ar^{f \tensfun{A} \coring{C}}[rr] & & N
\tensdia{A} \coring{C}}
\end{equation}
The $k$--module of all right $\coring{C}$--comodule morphisms from
$M$ to $N$ is denoted by $\hom{\coring{C}}{M}{N}$. The $k$--linear
category of all right $\coring{C}$--comodules will be denoted by
$\rcomod{\coring{C}}$. When $\coring{C} = A$, with the trivial
coring structure, the category $\rcomod{A}$ is just the category
of all right $A$--modules, which is `traditionally' denoted by
$\rmod{A}$.

Coproducts and cokernels in $\rcomod{\coring{C}}$ do exist and can
be already computed in $\rmod{A}$. Therefore,
$\rcomod{\coring{C}}$ has arbitrary inductive limits. If
${}_A\coring{C}$ is a flat module, then $\rcomod{\coring{C}}$ is
easily proved to be an abelian category.
\end{punto}

\begin{punto}\textbf{Bicomodules.}
Let $\rho_M : M \rightarrow M \tensor{A} \coring{C}$ be a comodule
structure over an $A'-A$--bimodule $M$, and assume that $\rho_M$
is $A'$--linear. For any right $A'$--module $X$, the right
$A$--linear map $X \tensfun{A'} \rho_{M} : X \tensor{A'} M
\rightarrow X \tensor{A'} M \tensor{A} \coring{C}$ makes $X
\tensor{A'} M$ a right $\coring{C}$--comodule. This leads to an
additive functor $- \tensor{A'} M : \rmod{A'} \rightarrow
\rcomod{\coring{C}}$. When $A' = A$ and $M = \coring{C}$, the
functor $- \tensor{A} \coring{C}$ is right adjoint to the
forgetful functor $U_A : \rcomod{\coring{C}} \rightarrow \rmod{A}$
(see \cite[Proposition 3.1]{Guzman:1989}, \cite[Lemma
3.1]{Brzezinski:2000unp}). Now assume that the $A'-A$--bimodule
$M$ is also a left $\coring{C}'$--comodule with structure map
$\lambda_M : M \rightarrow \coring{C}' \tensor{A} M$. It is clear
that $\rho_M : M \rightarrow M \tensor{A} \coring{C}$ is a
morphism of left $\coring{C}'$--comodules if and only if
$\lambda_M : M \rightarrow \coring{C}' \tensor{A'} M$ is a
morphism of right $\coring{C}$--comodules. In this case, we say
that $M$ is a $\coring{C}'-\coring{C}$--bicomodule. The
$\coring{C}'-\coring{C}$--bicomodules are the objects of a
$k$--linear category $\bcomod{\coring{C}'}{\coring{C}}$ whose
morphisms are those $A'-A$--bimodule maps which are morphisms of
$\coring{C}'$--comodules and of $\coring{C}$--comodules. Some
particular cases are now of interest. For instance, when
$\coring{C}' = A'$, the objects of the category
$\bcomod{A'}{\coring{C}}$ are the $A'-A$--bimodules with a right
$\coring{C}$--comodule structure $\rho_M : M \rightarrow M
\tensor{A} \coring{C}$ which is $A'$--linear.
\end{punto}

We will use the following lemma.

\begin{lemma}
If $M$ and $N$ are $A$-bimodules, and $P$ and $L$ are
$B$-bimodules, then
\[
(M\tensor{A}N)\tensor{k}(L\tensor{B}P) \cong
(M\tensor{k}L)\tensor{A\tensor{k}B}(N\tensor{k}P),
\]
and this is a natural isomorphism of $A \tensor{k} B$--bimodules.
\end{lemma}
\begin{proof}
Straightforward.
\end{proof}

With the aid of this lemma, we can give the tensor product of
corings as follows.

\begin{proposition}\label{tensorcoring}
Let $(\coring{C},\Delta_{\coring{C}}, \epsilon_{\coring{C}})$ and
$(\coring{D},\Delta_{\coring{D}}, \epsilon_{\coring{D}})$ be
corings over the algebras $A$ and $B$, respectively. Then
$\coring{C} \tensor{k} \coring{D}$ is an $A \tensor{k}
B$--coring, with the comultiplication
\[
\xymatrix{\coring{C}\otimes_{k}\coring{D}
\ar^(.25){\Delta_{\coring{C}} \otimes_{k} \Delta_{\coring{D}}}[rr]
 & &
 (\coring{C}\otimes_{A}\coring{C})\otimes_{k}(\coring{D}\otimes_{B}\coring{D})
  \cong (\coring{C}\otimes_{k}\coring{D})\otimes_{A\otimes_{k}B}(\coring{C}\otimes_{k}\coring{D})},
 \]
and the counit
\[
\xymatrix{\coring{C} \tensor{k} \coring{D}
\ar^{\epsilon_{\coring{C}} \tensor{k} \epsilon_{\coring{D}}}[rr] &
& A \tensor{k} B}
\]
\end{proposition}
\begin{proof}
Clearly, the proposed comultiplication and counit are
homomorphisms of $A \tensor{k} B$--bimodules. To check the
coassociative property, consider the diagram
\[
\xymatrix@C=-5pt{ &%
(\coring{C}\tensfun{A}\coring{C})
\tensfun{}(\coring{D}\tensfun{B}\coring{D})\ar[dr]^-\cong
\ar[ddd]|{(\coring{C}\tensfun{A}
\Delta)\tensfun{}(\coring{D}\tensfun{B}\Delta)} &%
\\
\coring{C}\tensfun{}\coring{D}\ar[ur]^-{\Delta\tensfun{}
\Delta}\ar[d]|{\Delta\tensfun{}\Delta} &%
 &
(\coring{C}\tensfun{}\coring{D}) \tensfun{A\tensfun{}B}
(\coring{C}\tensfun{}\coring{D})\ar[d]|{(\coring{C}\tensfun{}\coring{D})
\tensfun{A\tensfun{}B}(\Delta\tensfun{}\Delta)}\\%
(\coring{C}\tensfun{A}\coring{C})\tensfun{}
(\coring{D}\tensfun{B}\coring{D})\ar[dd]^\cong\ar[dr]^-
{(\Delta\tensfun{A}\coring{C})\tensfun{}(\Delta\tensfun{B}\coring{D})}
& & (\coring{C}\tensfun{}\coring{D})
\tensfun{A\tensfun{}B}(\coring{C}\tensfun{A}\coring{C})\tensfun{}
(\coring{D}\tensfun{B}\coring{D})\ar[dd]^\cong\\
 & (\coring{C}\tensfun{A}\coring{C}\tensfun{A}\coring{C})\tensfun{}(\coring{D}\tensfun{B}\coring{D}
\tensfun{B}\coring{D})\ar[ur]^-\cong\ar[dd]^\cong & \\
(\coring{C}\tensfun{}\coring{D})
\tensfun{A\tensfun{}B}(\coring{C}\tensfun{}\coring{D})\ar[dr]^-{(\Delta
\tensfun{}\Delta)\tensfun{A\tensfun{}B}(\coring{C}\tensfun{}\coring{D})}&
&
(\coring{C}\tensfun{}\coring{D})\tensfun{A\tensfun{}B}(\coring{C}\tensfun{}\coring{D})
\tensfun{A\tensfun{}B}(\coring{C}\tensfun{}\coring{D})\\
&
 (\coring{C}\tensfun{A}\coring{C})\tensfun{}(\coring{D}\tensfun{B}\coring{D})
\tensfun{A\tensfun{}B}(\coring{C}\tensfun{}\coring{D})\ar[ur]_-\cong
& & ,}
\]
whose commutativity follows from that of the inner diagrams. The
counitary property follows analogously from the commutativity of
the diagram
\[
\xymatrix{\coring{C}\tensfun{}\coring{D}\ar[r]^-{\Delta\tensfun{}
\Delta}\ar[dr]_{\cong} & (\coring{C}\tensfun{A}\coring{C})
\tensfun{}(\coring{D}\tensfun{B}\coring{D})\ar[r]^-\cong\ar[d]|{(\coring{C}\tensfun{A}
\epsilon)\tensfun{}(\coring{D}\tensfun{B}\epsilon)} &
(\coring{C}\tensfun{}\coring{D}) \tensfun{A\tensfun{}B}
(\coring{C}\tensfun{}\coring{D})\ar[d]|{(\coring{C}\tensfun{}\coring{D})
\tensfun{A\tensfun{}B}(\epsilon\tensfun{}\epsilon)}\\
 & (\coring{C}\tensfun{A}A)\tensfun{}(\coring{D}\tensfun{B}A)\ar[r]^-
{\cong} &
 (\coring{C}\tensfun{}\coring{D})
\tensfun{A\tensfun{}B}(A\tensfun{}B)}
\]
\end{proof}

\begin{example}\label{tensorentwining}
One important source of examples of corings are the entwining
structures, introduced by T. Brzezi\'nski and S. Majid in
\cite{Brzezinski/Majid:1998}. Let $A$, $C$ be a $k$--algebra and a
$k$--coalgebra, respectively. Assume that $A \tensor{} C$ is
endowed with an $A$--bimodule structure, where the left
$A$--module structure is the canonical one. The right $A$--module
structure induces a $k$--linear map $\psi : C \tensor{} A
\rightarrow A \tensor{} C$ given by $\psi(c \otimes a) = (1
\otimes c)a$.  Consider the maps
\[
\xymatrix{ A \tensor{} C \ar^-{A \tensor{} \Delta_C}[rr] & & A
\tensor{} C \tensor{} C \cong (A \tensor{} C) \tensor{A} A
\tensor{} C },
\]
and
\[
\xymatrix{A \tensor{} C \ar^-{A \epsilon_C}[rr] & & A \tensor{}k
\cong A }.
\]
By \cite[Proposition 2.2]{Brzezinski:2000unp} they give a
structure of $A$--coring on $A \tensor{} C$ if and only if
$(A,C,\psi)$ is an \emph{entwining structure} over $k$. The tensor
product of corings given in Proposition \ref{tensorcoring}
provides a canonical definition of tensor product of entwining
structures as follows. Let $(B,D,\varphi)$ be a second entwining
structure over $k$, and consider the tensor product $A \tensor{}
B$--coring $A \tensor{} C \tensor{} B \tensor{} D$. The
isomorphism of $k$--modules
\begin{equation}\label{isoc}
A \tensor{} C \tensor{} B \tensor{} D \cong A \tensor{} B
\tensor{} C \tensor{} D
\end{equation}
allows to induce an $A \tensor{} B$--bimodule structure on $A
\tensor{} B \tensor{} C \tensor{} D$ from the bimodule structure
of the $A \tensor{} B$--coring $A \tensor{} C \tensor{} B
\tensor{} D$. The resulting left $A \tensor{} B$--module structure
is then the canonical one, while the right $A \tensor{} B$
structure is, after some computations, the induced by the map
\[
\xymatrix{C \tensor{} D \tensor{} A \tensor{} B \cong C \tensor{}
A \tensor{} D \tensor{} B \ar^-{\psi \tensor{} \varphi}[rr] & & A
\tensor{} C \tensor{} B \tensor{} D \cong A \tensor{} B \tensor{}
C \tensor{} D },
\]
which will be denoted by $\psi \tensor{} \varphi$. It can be
checked in a straightforward way that the comultiplication and the
counit of the tensor product $A \tensor{} B$--coring $A \tensor{}
C \tensor{} B \tensor{} D$ define, by using the isomorphism
\eqref{isoc}, a comultiplication and a counit on $A \tensor{} B
\tensor{} C \tensor{} D$ in such a way that $(A \tensor{} B, C
\tensor{} D, \psi \tensor{} \varphi)$ is an entwining structure,
which we call \emph{tensor product entwining structure}. Of
course, \eqref{isoc} becomes an isomorpism of $A \tensor{}
B$--corings.
\end{example}

\begin{punto}\textbf{Opposite coring.}
The notation $A^{\circ}$ stands for the opposite algebra of $A$;
every $A-B$--bimodule $M$ gives in an obvious way a
$B^{\circ}-A^{\circ}$ bimodule, denoted by $M^{\circ}$. So every
$A$--coring $\coring{C}$ defines an $A^{\circ}$--coring
$\coring{C}^{\circ}$ as follows: we have an isomorphism of
$k$--modules $\tau : \coring{C} \tensor{A} \coring{C} \rightarrow
\coring{C}^{\circ} \tensor{A^{\circ}} \coring{C}^{\circ}$ defined
by $\tau (x \tensor{A} y) = y \tensor{A^{\circ}} x $. Then the
$k$--linear map
\[
\xymatrix{\coring{C} \ar^(.4){\Delta}[r] & \coring{C} \tensor{A}
\coring{C} \ar^(.4){\tau}[r] & \coring{C}^{\circ}
\tensor{A^{\circ}} \coring{C}^{\circ}}
\]
can be viewed as an $A^{\circ}-A^{\circ}$--bimodule map
$\Delta^{\circ} : \coring{C}^{\circ} \rightarrow
\coring{C}^{\circ} \tensor{A^{\circ}} \coring{C}^{\circ}$, which
is clearly coassociative. Consider the counit $\epsilon :
\coring{C} \rightarrow A$ as an $A^{\circ}-A^{\circ}$--bimodule
map $\epsilon^{\circ} : \coring{C}^{\circ} \rightarrow A^{\circ}$
to obtain an $A^{\circ}$--coring
$(\coring{C}^{\circ},\Delta^{\circ},\epsilon^{\circ})$ called the
\emph{opposite coring} of $\coring{C}$.
\end{punto}

We are now in position to show that bicomodules are in fact
comodules.

\begin{proposition}\label{comodbicomod}
Let $A, B$ be algebras, and $\coring{C}, \coring{D}$ be corings
over $A$ and $B$, respectively. Then there is an isomorphism of
categories $\bcomod{\coring{C}}{\coring{D}} \cong
\lcomod{\coring{C} \tensor{k} \coring{D}^{\circ}}$.
\end{proposition}
\begin{proof}
If $(M, \lambda_{M}, \rho_{M})$ is a
 $\coring{C}-\coring{D}^{\circ}$-bicomodule, then the
 homomorphism of $A\tensor{k}B$-bimodules
\[
\xymatrix{\lambda: M\ar[r]^-{\lambda_{M}} & \coring{C}\tensor{A}M
\ar[r]^-{\coring{C}\tensor{A}\rho_{M}}
 &
 \coring{C}\tensor{A}M\tensor{B^{\circ}}\coring{D^{\circ}}\ar[r]^-{\cong}
 & (\coring{C}\tensor{k}\coring{D})\tensor{A\tensor{k}B}M}
 \]
 defines a structure of left
 $\coring{C}\tensor{k}\coring{D}$-comodule. The pseudo-coassociative property for
 $\lambda$ follows from the commutativity of all the inner diagrams
 in the following one:
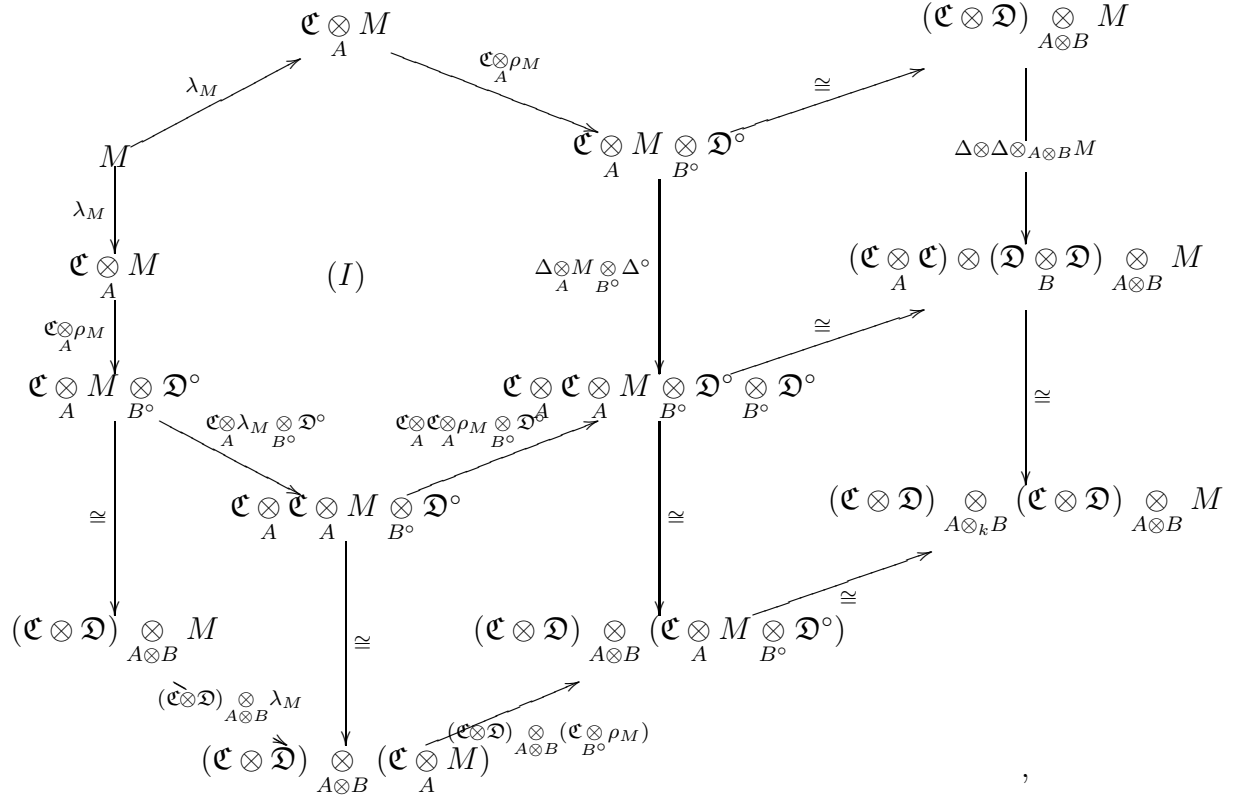
\[
\xymatrix@C=-7pt@W=0pt@M=0pt{ & \coring{C}\tensfun{A}M
\ar[dr]^-{\coring{C}\tensfun{A}\rho_{M}}&
&(\coring{C}\tensfun{}\coring{D})\tensfun{A\tensfun{}B}M
 \ar[dd]|-{\Delta\tensfun{}\Delta\tensor{A\tensfun{}B}M} \\
          M\ar[ur]^-{\lambda_{M}}\ar[d]_-{\lambda_{M}} & &
          \coring{C}\tensfun{A}M\tensfun{B^{\circ}}\coring{D^{\circ}}\ar[ur]^-{\cong}
 \ar[dd]_-{\Delta\tensfun{A}M\tensfun{B^{\circ}}\Delta^{\circ}}& \\
           \coring{C}\tensfun{A}M\ar[d]_-{\coring{C}\tensfun{A}\rho_{M}} & (I) & &
           (\coring{C}\tensfun{A}\coring{C})\tensfun{}(\coring{D}\tensfun{B}
\coring{D})\tensfun{A\tensfun{}B}M\ar[dd]^-{\cong}\\
           \coring{C}\tensfun{A}M\tensfun{B^\circ}\coring{D}^{\circ}
  \ar[dr]^-{\coring{C}\tensfun{A}\lambda_{M}\tensfun{B^\circ}\coring{D}^{\circ}}
  \ar[dd]_-{\cong}& & \coring{C}\tensfun{A}\coring{C}\tensfun{A}M\tensfun{B^{\circ}}
 \coring{D^{\circ}}\tensfun{B^{\circ}}\coring{D}^{\circ}\ar[ur]^-{\cong}
 \ar[dd]^-{\cong}& \\
           & \coring{C}\tensfun{A}\coring{C}\tensfun{A}M\tensfun{B^{\circ}}
 \coring{D^{\circ}}
 \ar[ur]^-{\coring{C}\tensfun{A}\coring{C}\tensfun{A}\rho_{M}
 \tensfun{B^{\circ}}\coring{D}^{\circ}}
 \ar[dd]^-{\cong}& &(\coring{C}\tensfun{}\coring{D})\tensfun{A\tensor{k}B}
 (\coring{C}\tensfun{}\coring{D})\tensfun{A\tensfun{}B}M \\
  (\coring{C}\tensfun{}\coring{D})\tensfun{A\tensfun{}B}M
 \ar[dr]|-{(\coring{C}\tensfun{}\coring{D})\tensfun{A\tensfun{}B}\lambda_{M}}
 & & (\coring{C}\tensfun{}\coring{D})\tensfun{A\tensfun{}B}
 (\coring{C}\tensfun{A}M\tensfun{B^{\circ}}\coring{D^{\circ}})\ar[ur]_-{\cong}& \\
           & (\coring{C}\tensfun{}\coring{D})\tensfun{A\tensfun{}B}(\coring{C}\tensfun{A}M)
\ar[ur]_-{(\coring{C}\tensfun{}\coring{D})\tensfun{A\tensfun{}B}
(\coring{C}\tensfun{B^{\circ}}\rho_{M})}& &  ,}
 \]
Let us check that $(I)$ commutes, which perhaps is not evident:

\begin{eqnarray*}
(\Delta \tensfun{A} M \tensfun{B^{\circ}}
\Delta^{\circ})(\coring{C} \tensfun{A} \rho_M)\lambda_M & = &
(\Delta \tensfun{A} \rho_M \tensfun{B^{\circ}}
\coring{D}^{\circ})(\coring{C} \tensfun{A}
\rho_M)\lambda_M \\
 & = & (\Delta \tensfun{A} \rho_M \tensfun{B^{\circ}}
 \coring{D}^{\circ})(\lambda_M \tensfun{B^{\circ}}
 \coring{D})\rho_M \\
  & = & (\coring{C} \tensfun{A} \coring{C} \tensor{A} \rho_M
  \tensfun{B^{\circ}} \coring{D}^{\circ})
  (\Delta \tensfun{A} M \tensfun{B^{\circ}} \coring{D}^{\circ})
  (\lambda_M \tensfun{B^{\circ}} \coring{D})\rho_M \\
  & = & (\coring{C} \tensfun{A} \coring{C} \tensfun{A} \rho_M
  \tensfun{B^{\circ}} \coring{D}^{\circ})(\coring{C} \tensfun{A}
  \lambda_M \tensfun{B^{\circ}} \coring{D}^{\circ})(\coring{C}
  \tensfun{A} \rho_M)\lambda_M
\end{eqnarray*}

Since the diagram
 \[
\xymatrix{M\ar[r]^-{\lambda_{M}}\ar[dr]_-{\cong} &
\coring{C}\tensfun{A}M
\ar[r]^-{\coring{C}\tensfun{A}\rho_{M}}\ar[d]^-{\epsilon\tensfun{A}M}
 &
 \coring{C}\tensfun{A}M\tensor{B^{\circ}}\coring{D^{\circ}}\ar[r]^-{\cong}
 \ar[d]^-{\epsilon\tensfun{A}M\tensfun{B^{\circ}}\epsilon^{\circ}}
 &
 (\coring{C}\tensor{k}\coring{D})\tensfun{A\tensfun{}B}M
 \ar[d]^-{(\epsilon\tensfun{}\epsilon)\tensfun{A\tensfun{}B}M}
 \\ & A\tensor{A}M
\ar[r]^-{\cong}
 &
 A\tensfun{A}M\tensfun{B^{\circ}}B^{\circ}\ar[r]^-{\cong}
 & (A\tensfun{}B)\tensfun{A\tensfun{}B}M}
 \]
 is also commutative, the counital property for the coaction $\lambda$ holds, and we have that
 $(M,\lambda)$ is a left $\coring{C}$--comodule.
 Moreover, it is easy to see that every homomorphism of
 $\coring{C}-\coring{D}^{\circ}$--bicomodules becomes a morphism
 of left $\coring{C} \tensor{k} \coring{D}$--comodules. \\
 Conversely, let $(M,\lambda)$ any left $\coring{C} \tensor{k}
 \coring{D}$--comodule. Then the commutativity of the diagrams
 \[
\xymatrix@C=-7pt@M=0pt{M\ar[dr]^-{\lambda}\ar[dd]_-{\lambda}& &
\coring{C}\tensfun{A}M\tensfun{B^{\circ}}\coring{D}^{\circ}
  \ar[dr]^-{\coring{C}\tensfun{A}M\tensfun{B^{\circ}}\epsilon^{\circ}}
  \ar[dd]_-{\Delta\tensfun{A}M\tensfun{B^{\circ}}\Delta^{\circ}}& \\
          & (\coring{C}\tensfun{}\coring{D})\tensfun{A\tensfun{}B}M\ar[ur]^-{\cong}
 \ar[dd]^-{\Delta\tensfun{A\tensfun{}B}M}& & \coring{C}\tensfun{A}M\ar[dddd]^-{\Delta\tensfun{A}M}\\
         (\coring{C}\tensfun{}\coring{D})\tensfun{A\tensfun{}B}M
 \ar[dr]^-{(\coring{C}\tensfun{}\coring{D})\tensfun{A\tensfun{}B}\lambda}
 \ar[dd]_-{\cong} &
 &\coring{C}\tensfun{A}\coring{C}\tensfun{A}M\tensfun{B^{\circ}}\coring{D}^{\circ}
 \tensfun{B^{\circ}}\coring{D}^{\circ}
 \ar[rddd]|-{ \coring{C}\tensfun{A}\coring{C}
 \tensfun{A}M\tensfun{B^{\circ}}\epsilon^{\circ}
 \tensfun{B^{\circ}}\epsilon^{\circ}}
 \ar[dddd]_-{
 \coring{C}\tensfun{A}\coring{C}\tensfun{A}M\tensfun{B^{\circ}}\coring{D}^{\circ}
 \tensfun{B^{\circ}}\epsilon^{\circ}} & \\
          & (\coring{C}\tensfun{}\coring{D})\tensfun{A\tensfun{}B}
 (\coring{C}\tensfun{}\coring{D})\tensfun{A\tensfun{}B}M
\ar[ur]^-{\cong}\ar[dd]_-{\cong}& & \\
          \coring{C}\tensfun{A}M\tensfun{B^{\circ}}\coring{D}^{\circ}
 \ar[dr]^-{\coring{C}\tensfun{A}\lambda\tensfun{B^{\circ}}\coring{D}^{\circ}}
 \ar[dd]_-{\coring{C}\tensfun{A}M\tensfun{B^{\circ}}\epsilon^{\circ}}& & & \\
          &\coring{C}\tensfun{A}((\coring{C}\tensfun{}\coring{D})
  \tensfun{A\tensfun{}B}M)\tensfun{B^{\circ}}\coring{D}^{\circ}
  \ar[dd]|-{ \coring{C}\tensfun{A}((\coring{C}\tensfun{}\coring{D})
  \tensfun{A\tensfun{}B}M)\tensfun{B^{\circ}}\epsilon^{\circ}} & & \coring{C}\tensfun{A}\coring{C}\tensfun{A}M\\
         \coring{C}\tensfun{A}M\ar[dr]_-{\coring{C}\tensfun{A}\lambda} & & (\coring{C}\tensfun{A}\coring{C})\tensfun{A}M
   \tensfun{B^{\circ}}\coring{D}^{\circ}
   \ar[ur]|-{(\coring{C}\tensfun{A}\coring{C})\tensfun{A}M
   \tensfun{B^{\circ}}\epsilon^{\circ}}& \\
          & \coring{C}\tensfun{A}((\coring{C}\tensfun{}\coring{D})
  \tensfun{A\tensfun{}B}M)\ar[ur]^-{\cong}& &, }
    \]
    and
    \[
\xymatrix{M\ar[r]^-{\lambda}\ar[dr]_-{\cong}
 &
 (\coring{C}\tensfun{}\coring{D})\tensfun{A\tensfun{}B}M\ar[r]^-{\cong}
 \ar[d]^-{(\epsilon\tensfun{}\epsilon)\tensfun{A\tensfun{}B}M}
  &
  \coring{C}\tensfun{A}M\tensfun{B^{\circ}}\coring{D}^{\circ}
  \ar[r]^-{\coring{C}\tensfun{A}M\tensfun{B^{\circ}}\epsilon^{\circ}}
  \ar[d]^-{\epsilon\tensfun{A}M\tensfun{B^{\circ}}\epsilon^{\circ}}
   &
   \coring{C}\tensfun{A}M\ar[d]^-{\epsilon\tensfun{A}M}
    \\ & (A\tensfun{}B)\tensfun{A\tensfun{}B}M\ar[r]^-{\cong}
   &  A\tensfun{A}M\tensfun{B^{\circ}}B^{\circ}
  \ar[r]^-{\cong}
   &
   A\tensfun{A}M}
   \]
   show that $(M,\lambda_M)$ is a left $\coring{C}$--comodule with
   structure map
   \[
   \lambda_{M}: \xymatrix{M\ar[r]^-{\lambda}
 &
 (\coring{C}\tensor{k}\coring{D})\tensor{A\tensor{k}B}M\ar[r]^-{\cong}
  &
  \coring{C}\tensor{A}M\tensor{B^{\circ}}\coring{D}^{\circ}
  \ar[rr]^-{\coring{C}\tensor{A}M\tensor{B^{\circ}}\epsilon^{\circ}}
  & &  \coring{C}\tensor{A}M}.
  \]
Analogously, the map
\[
\rho_{M}: \xymatrix{M\ar[r]^-{\lambda}
 &
 (\coring{C}\tensor{k}\coring{D})\tensor{A\tensor{k}B}M\ar[r]^-{\cong}
  &
  \coring{C}\tensor{A}M\tensor{B^{\circ}}\coring{D}^{\circ}
  \ar[rr]^-{\epsilon\tensor{A}M\tensor{B^{\circ}}\coring{D}^{\circ}}
  & &  M\tensor{B^{\circ}}\coring{D}^{\circ}}
  \]
  endows $M$ with a structure of a right
  $\coring{D}^{\circ}$--comodule. Next, we shall check that
  $(M,\lambda_M,\rho_M)$ is a
  $\coring{C}-\coring{D}^{\circ}$--bicomodule. Consider for this
  the commutative diagrams
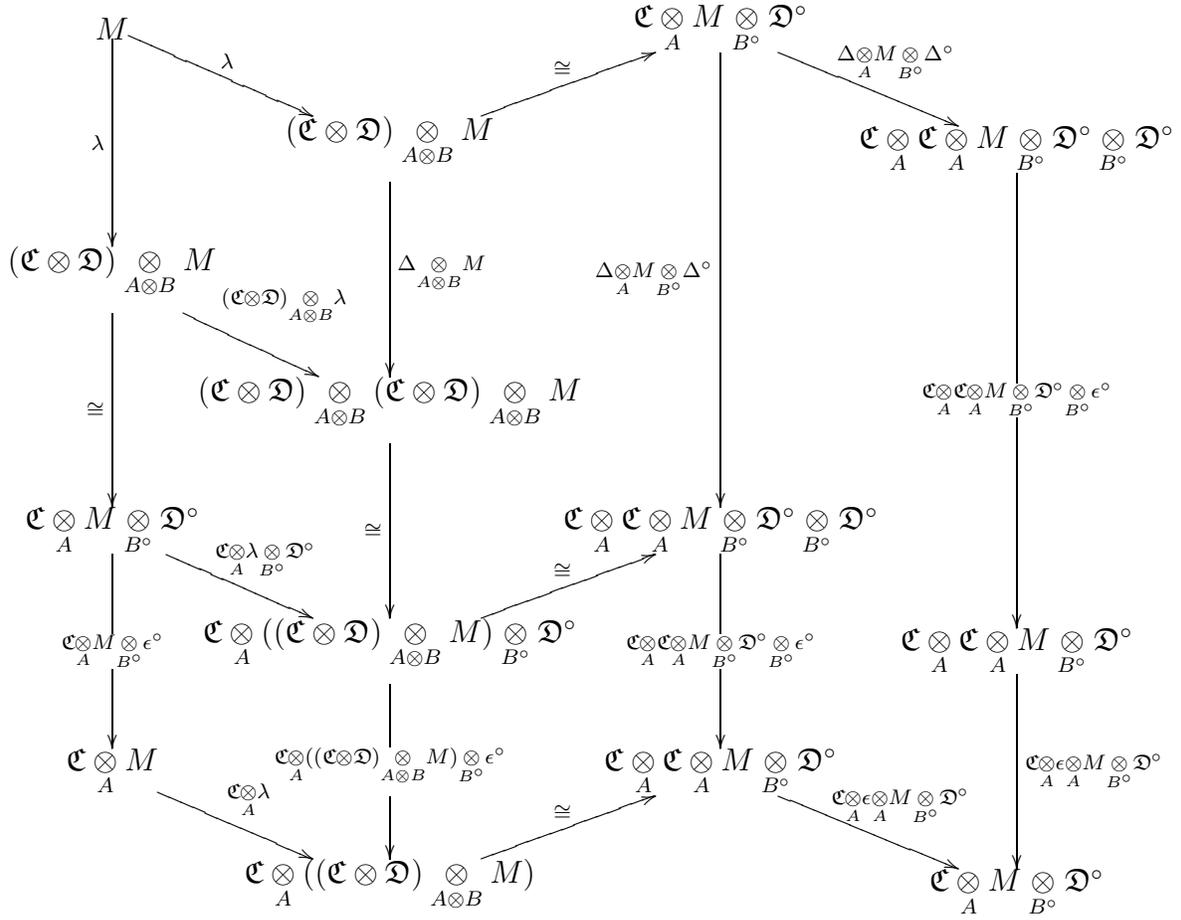
\[
  \xymatrix@C=-7pt@M=0pt{M\ar[dr]^-{\lambda}\ar[dd]_-{\lambda} &
    &\coring{C}\tensfun{A}M\tensfun{B^{\circ}}\coring{D}^{\circ}
  \ar[dr]^-{\Delta\tensfun{A}M\tensfun{B^{\circ}}\Delta^{\circ}}
  \ar[dddd]_-{\Delta\tensfun{A}M\tensfun{B^{\circ}}\Delta^{\circ}} & \\
                      & (\coring{C}\tensfun{}\coring{D})\tensfun{A\tensfun{}B}M\ar[ur]^-{\cong}
 \ar[dd]^-{\Delta\tensfun{A\tensfun{}B}M}& &
 \coring{C}\tensfun{A}\coring{C}\tensfun{A}M\tensfun{B^{\circ}}\coring{D}^{\circ}
 \tensfun{B^{\circ}}\coring{D}^{\circ}\ar[dddd]|-{\coring{C}\tensfun{A}\coring{C}\tensfun{A}M\tensfun{B^{\circ}}\coring{D}^{\circ}
 \tensfun{B^{\circ}}\epsilon^{\circ}}\\
                      (\coring{C}\tensfun{}\coring{D})\tensfun{A\tensfun{}B}M
 \ar[dr]^-{(\coring{C}\tensfun{}\coring{D})\tensfun{A\tensfun{}B}\lambda}
 \ar[dd]_-{\cong}& & & \\
                      & (\coring{C}\tensfun{}\coring{D})\tensfun{A\tensfun{}B}
 (\coring{C}\tensfun{}\coring{D})\tensfun{A\tensfun{}B}M
 \ar[dd]_-{\cong} & & \\
                     \coring{C}\tensfun{A}M\tensfun{B^{\circ}}\coring{D}^{\circ}
 \ar[dr]^-{\coring{C}\tensfun{A}\lambda\tensfun{B^{\circ}}\coring{D}^{\circ}}
 \ar[dd]|-{\coring{C}\tensfun{A}M\tensfun{B^{\circ}}\epsilon^{\circ}} &
 &\coring{C}\tensfun{A}\coring{C}\tensfun{A}M\tensfun{B^{\circ}}\coring{D}^{\circ}
 \tensfun{B^{\circ}}\coring{D}^{\circ}
 \ar[dd]|-{\coring{C}\tensfun{A}\coring{C}\tensfun{A}M\tensfun{B^{\circ}}\coring{D}^{\circ}
 \tensfun{B^{\circ}}\epsilon^{\circ}} & \\
                      & \coring{C}\tensfun{A}((\coring{C}\tensfun{}\coring{D})
  \tensfun{A\tensfun{}B}M)\tensfun{B^{\circ}}\coring{D}^{\circ}
  \ar[ur]^-{\cong}
  \ar[dd]|-{ \coring{C}\tensfun{A}((\coring{C}\tensfun{}\coring{D})
  \tensfun{A\tensfun{}B}M)\tensfun{B^{\circ}}\epsilon^{\circ}}& &
  \coring{C}\tensfun{A}\coring{C}\tensfun{A}M\tensfun{B^{\circ}}\coring{D}^{\circ}
   \ar[dd]^-{\coring{C}\tensfun{A}\epsilon\tensfun{A}M\tensfun{B^{\circ}}\coring{D}^{\circ}}\\
                     \coring{C}\tensfun{A}M\ar[dr]^-{\coring{C}\tensfun{A}\lambda} & &\coring{C}\tensfun{A}\coring{C}\tensfun{A}M
   \tensfun{B^{\circ}}\coring{D}^{\circ}
   \ar[dr]^-{\coring{C} \tensfun{A}\epsilon\tensfun{A}M
   \tensfun{B^{\circ}}\coring{D}^{\circ}} & \\
                      & \coring{C}\tensfun{A}((\coring{C}\tensfun{}\coring{D})
  \tensfun{A\tensfun{}B}M)\ar[ur]^-{\cong}& & \coring{C}\tensfun{A}M\tensfun{B^{\circ}}\coring{D}^{\circ}}
  \]
  and
  \[
\xymatrix@C=-7pt@M=0pt{ M\ar[dr]^-{\lambda}\ar[dd]_-{\lambda}& &
\coring{C}\tensfun{A}M\tensfun{B^{\circ}}\coring{D}^{\circ}
  \ar[dr]^-{\Delta\tensfun{A}M\tensfun{B^{\circ}}\Delta^{\circ}}
  \ar[dddd]_-{\Delta\tensfun{A}M\tensfun{B^{\circ}}\Delta^{\circ}}& \\
                  &(\coring{C}\tensfun{}\coring{D})\tensfun{A\tensfun{}B}M\ar[ur]^-{\cong}
 \ar[dd]^-{\Delta\tensfun{A\tensfun{}B}M} &
 &\coring{C}\tensfun{A}\coring{C}\tensfun{A}M\tensfun{B^{\circ}}\coring{D}^{\circ}
 \tensfun{B^{\circ}}\coring{D}^{\circ}\ar[dddd]|-{\coring{C}\tensfun{A}\coring{C}\tensfun{A}M\tensfun{B^{\circ}}\coring{D}^{\circ}
 \tensfun{B^{\circ}}\epsilon^{\circ}} \\
                  (\coring{C}\tensfun{}\coring{D})\tensfun{A\tensfun{}B}M
 \ar[dr]^-{(\coring{C}\tensfun{}\coring{D})\tensfun{A\tensfun{}B}\lambda}
 \ar[dd]_-{\cong}& & & \\
                  & (\coring{C}\tensfun{}\coring{D})\tensfun{A\tensfun{}B}
 (\coring{C}\tensfun{}\coring{D})\tensfun{A\tensfun{}B}M
 \ar[dd]_-{\cong}& & \\
                 \coring{C}\tensfun{A}M\tensfun{B^{\circ}}\coring{D}^{\circ}
 \ar[dr]^-{\coring{C}\tensfun{A}\lambda\tensfun{B^{\circ}}\coring{D}^{\circ}}
 \ar[dd]_-{\epsilon\tensfun{A}M\tensfun{B^{\circ}}\coring{D}^{\circ}} & &
 \coring{C}\tensfun{A}\coring{C}\tensfun{A}M\tensfun{B^{\circ}}\coring{D}^{\circ}
 \tensfun{B^{\circ}}\coring{D}^{\circ}
 \ar[dd]|-{\epsilon\tensfun{A}\coring{C}\tensfun{A}M\tensfun{B^{\circ}}\coring{D}^{\circ}
 \tensfun{B^{\circ}}\coring{D}^{\circ}} & \\
                  & \coring{C}\tensfun{A}((\coring{C}\tensfun{}\coring{D})
  \tensfun{A\tensfun{}B}M)\tensfun{B^{\circ}}\coring{D}^{\circ}
  \ar[ur]^-{\cong}
  \ar[dd]|-{ \epsilon\tensfun{A}((\coring{C}\tensfun{}\coring{D})
  \tensfun{A\tensfun{}B}M)\tensfun{B^{\circ}}\coring{D}^{\circ}}& &
  \coring{C}\tensfun{A}\coring{C}\tensfun{A}M\tensfun{B^{\circ}}\coring{D}^{\circ}
   \ar[dd]^-{\coring{C}\tensfun{A}\epsilon\tensfun{A}M\tensfun{B^{\circ}}\coring{D}^{\circ}}\\
                 M\tensfun{B^{\circ}}\coring{D}^{\circ}
   \ar[dr]^-{\lambda\tensfun{B^{\circ}}\coring{D}^{\circ}} & &\coring{C}\tensfun{A}M
   \tensfun{B^{\circ}}\coring{D}^{\circ}\tensfun{B^{\circ}}\coring{D}^{\circ}
   \ar[dr]^-{\coring{C}\tensfun{A}M
   \tensfun{B^{\circ}}\epsilon^{\circ}\tensfun{B^{\circ}}\coring{D}^{\circ}} & \\
                  & ((\coring{C}\tensfun{}\coring{D})
  \tensfun{A\tensfun{}B}M) \tensfun{B^{\circ}} \coring{D}^{\circ} \ar[ur]^-{\cong}&
   &\coring{C}\tensfun{A}M\tensfun{B^{\circ}}\coring{D}^{\circ}}
    \]
    which say that $(\coring{C} \tensor{A} \rho_M)\lambda_M =
    (\lambda_M \tensor{B^{\circ}} \coring{D}^{\circ})\rho_M$, that is,
    $M$ is a bicomodule.  Moreover, it is easy to see that every
    homomorphism of left $\coring{C} \tensor{k}
    \coring{D}$--comodules is a homomorphism of $\coring{C} -
    \coring{D}^{\circ}$--bicomodules. We have constructed two
    functors which are easily seen to be mutually inverse.
    This gives the desired isomorphism of categories.
\end{proof}

\section{Coseparable corings and semisimple corings}

Recall from \cite{DeMeyer/Ingraham:1971} that the algebra $A$ over
the commutative ring $k$ is said to be \emph{separable} if the
multiplication map $\mu : A \tensor{k} A \rightarrow A$ splits as
an $A$--bimodule epimorphism, i.e., there is a homomorphism of
$A$--bimodules $\theta : A \rightarrow A \tensor{k} A$ such that
$\mu \theta = A$. Dually, one obtains the notion of a coseparable
coalgebra and, more generally, of a coseparable coring. Following
\cite{Guzman:1989}, we say that an $A$--coring $\coring{C}$ is
\emph{coseparable} if there is a $\coring{C}$--bicomodule map $\pi
: \coring{C} \tensor{A} \coring{C} \rightarrow \coring{C}$ such
that $\pi \Delta = \coring{C}$ (here, $A$ is not assumed to be
separable).

Both notions can be expressed functorially by means of the
separable functors introduced by C. N\u ast\u asescu, M. Van den
Bergh and F. Van Oystaeyen in
\cite{Nastasescu/VandenBergh/VanOystaeyen:1989}. Thus, $A$ is a
separable $k$--algebra if and only if the forgetful functor
$\rmod{A} \rightarrow \rmod{k}$ is separable \cite[Proposition
1.3]{Nastasescu/VandenBergh/VanOystaeyen:1989}, and $\coring{C}$
is a coseparable coring if and only if the forgetful functor
$\rcomod{C} \rightarrow \rmod{A}$ is separable \cite[Corollary 3.6
]{Brzezinski:2000unp}. We take advantage of this functorial
approach in our development (see \cite{Gomez:2001unp} for a
general treatment of separable functors for corings). Separable
functors are nicely characterized when they are in adjunction, see
\cite[Theorem 1.2]{Rafael:1990}. The first theorem in this section
needs a previous lemma.

\begin{lemma}
Let $\coring{C}$ and $\coring{D}$ be corings over $k$--algebras
$A$ and $B$, respectively. If $F :
\bcomod{\coring{C}}{\coring{D}} \rightarrow
\bcomod{A}{\coring{D}}$ and $F':\bcomod{A}{\coring{D}}
\rightarrow \rcomod{\coring{D}}$ are the forgetful functors, then
\begin{enumerate}[(1)]
\item The functor $G = \coring{C}  \tensor{A} - :
\bcomod{A}{\coring{D}} \rightarrow
\bcomod{\coring{C}}{\coring{D}}$ is right adjoint to the functor
$F$ with unit $\lambda : 1 \rightarrow GF $ given at every $M \in
\bcomod{\coring{C}}{\coring{D}}$ by the left
$\coring{C}$--comodule structure map $\lambda_M : M \rightarrow
\coring{C} \tensor{A} M$.
\item The functor $G' = A \tensor{k} - : \rcomod{\coring{D}}
\rightarrow \bcomod{A}{\coring{D}}$ is left adjoint to $F'$ with
the unit $c : G'F' \rightarrow 1$ defined at every $M \in
\bcomod{A}{\coring{D}}$ by the left $A$--module structure map
$c_M : A \tensor{k} M \rightarrow M$.
\end{enumerate}
\end{lemma}
\begin{proof}
(1) This is \cite[Proposition 3.1]{Guzman:1989}.\\
(2) Straightforward.
\end{proof}

\begin{theorem}\label{separaD}
Let $\coring{C}$ be a coseparable $A$--coring,
and assume that $A$ is a separable algebra over the commutative
ring $k$. If $\coring{D}$ is a coring over a $k$--algebra $B$,
then the forgetful functor $U : \bcomod{\coring{C}}{\coring{D}}
\rightarrow \rcomod{\coring{D}}$ is separable.
\end{theorem}
\begin{proof}
We first decompose the forgetful functor $U$ accordingly the
following diagram
\[
\xymatrix{\bcomod{\coring{C}}{\coring{D}} \ar^{U}[rr] \ar^{F}[dr]
& &
\rcomod{\coring{D}} \\
& \bcomod{A}{\coring{D}} \ar^{F'}[ur] & .}
\]
 We shall
prove that $F$ and $F'$ are separable functors, which implies, by
\cite[Lemma 1.1]{Nastasescu/VandenBergh/VanOystaeyen:1989}, that
$U$ is separable. By \cite[Theorem 1.2]{Rafael:1990}, $F$ is
separable if and only if there is a natural transformation $\mu :
GF \rightarrow 1$ such that $\mu \circ \lambda =
1_{\bcomod{\coring{C}}{\coring{D}}}$. By
\cite[3.6]{Brzezinski:2000unp}, there is a morphism of
$A$--bimodules $\gamma : \coring{C} \tensor{A} \coring{C}
\rightarrow A$ such that $\gamma \Delta = \epsilon$ and
$(\coring{C}\tensor{A}\gamma)(\Delta\tensor{A}\coring{C})=
 (\gamma\tensor{A}\coring{C})(\coring{C}\tensor{A}\Delta)$. Then, for
 every $\coring{C}-\coring{D}$--bicomodule
 $(M,\lambda_M,\rho_M)$, the commutativity of the following
 diagrams
 \[
 \xymatrix{\coring{C}\tensfun{A}M\ar[r]^-{\coring{C}\tensfun{A}\lambda_{M}}
 \ar[d]_-{\coring{C}\tensfun{A}\rho_{M}} &
 \coring{C}\tensfun{A}\coring{C}\tensfun{A}M\ar[r]^-{\gamma\tensfun{A}M}
 \ar[d]^-{\coring{C}\tensfun{A}\coring{C}\tensfun{A}\rho_{M}}
 & M\ar[d]^-{\rho_{M}}
  \\ \coring{C}\tensfun{A}M\tensfun{B}\coring{D}
  \ar[r]_-{\coring{C}\tensfun{A}\lambda_{M}\tensfun{B}\coring{D}}
   &
   \coring{C}\tensfun{A}\coring{C}\tensfun{A}M\tensfun{B}\coring{D}
   \ar[r]_-{\gamma\tensfun{A}M\tensfun{B}\coring{D}}
    & M\tensfun{B}\coring{D}}
   \]

 \[
 \xymatrix{\coring{C}\tensfun{A}M\ar[r]^-{\coring{C}\tensfun{A}\lambda_{M}}
 \ar[d]_-{\Delta\tensfun{A}M} &
 \coring{C}\tensfun{A}\coring{C}\tensfun{A}M\ar[r]^-{\gamma\tensfun{A}M}
  & M\ar[d]^-{\lambda_{M}}
  \\ \coring{C}\tensfun{A}\coring{C}\tensfun{A}M
  \ar[r]_-{\coring{C}\tensfun{A}\coring{C}\tensfun{A}\lambda_{M}}
   &
   \coring{C}\tensfun{A}\coring{C}\tensfun{A}\coring{C}\tensfun{A}M
   \ar[r]_-{\coring{C}\tensfun{A}\gamma\tensfun{A}M}
    & \coring{C}\tensfun{A}M}
   \]
   gives that $\mu_{M}: \xymatrix{\coring{C}\tensor{A}M\ar[r]^-{\coring{C}\tensor{A}
 \lambda_{M}} &
 \coring{C}\tensor{A}\coring{C}\tensor{A}M\ar[r]^-{\gamma\tensor{A}M}
  & M} $ is a homomorphism of
  $\coring{C}-\coring{D}$-bicomodules. On the other hand, if $f
  : M \rightarrow N$ is a homomorphism of
  $\coring{C}-\coring{D}$-bicomodules, then the diagram
  \[ \xymatrix{\coring{C}\tensfun{A}M\ar[r]^-{\coring{C}\tensfun{A}\lambda_{M}}
 \ar[d]_-{\coring{C}\tensfun{A}f} &
 \coring{C}\tensfun{A}\coring{C}\tensfun{A}M\ar[r]^-{\gamma\tensfun{A}M}
  & M\ar[d]^-{f}
  \\ \coring{C}\tensfun{A}N
  \ar[r]_-{\coring{C}\tensfun{A}\lambda_{N}}
   &
   \coring{C}\tensfun{A}\coring{C}\tensfun{A}N
   \ar[r]_-{\gamma\tensfun{A}N}
    & N}
   \]
   is commutative. Therefore we have defined a natural
   transformation
   $\mu:GF\rightarrow1_{\bcomod{\coring{C}}{\coring{D}}}$.
   Moreover, for every $M \in \bcomod{\coring{C}}{\coring{D}}$,
   one has $\mu_{M}\circ\lambda_{M}=(\gamma\tensor{A}M)(\Delta\tensor{A}M)
   \lambda_{M}=(\epsilon\tensor{A}M)\lambda_{M}=M$. Thus,
   $\lambda$ splits off, and $F$ is a separable functor.

   \noindent Now, to prove that $F'$ is a separable functor it
   suffices, by \cite[Theorem 1.2]{Rafael:1990}, to show that there is a natural
   transformation $\nu : 1 \rightarrow G'F'$ such that $c \circ
   \nu = 1$. Since $A$ is a separable $k$--algebra, we have an
   $A$--bimodule map $\theta : A \rightarrow A \tensor{k} A$ such
   that $\mu \theta = A$, where $\mu : A \tensor{k} A \rightarrow
   A$ is the multiplication map. Then for every $A-\coring{D}$-bicomodule
   $M$ we have the commutative diagram
   \[
     \xymatrix{M\ar[r]^-{\cong}\ar[d]_-{\rho_{M}} &
    A\tensfun{A}M\ar[r]^-{\theta\tensfun{A}M}\ar[d]^-{A\tensfun{A}\rho_{M}}
     & A\tensfun{}A\tensfun{A}M\ar[r]^-{\cong}
     \ar[d]^-{A\tensfun{}A\tensfun{A}\rho_{M}}
     & A\tensfun{}M\ar[d]^-{A\tensfun{}\rho_{M}}
      \\ M\tensfun{B}\coring{D}\ar[r]_-{\cong} &
      A\tensfun{A}M\tensfun{B}\coring{D}
     \ar[r]_-{\theta\tensfun{A}M\tensfun{B}\coring{D}}
       &
       A\tensfun{}A\tensfun{A}M\tensfun{B}\coring{D}\ar[r]_-{\cong}
        & A\tensfun{}M\tensfun{B}\coring{D},}
   \]
   which says that $\nu_{M}:\xymatrix{M\ar[r]^-{\cong} &
    A\tensor{A}M\ar[r]^-{\theta\tensor{A}M}
     & A\tensor{k}A\tensor{A}M\ar[r]^-{\cong}
     & A\tensor{k}M}$ is $\coring{D}$--colinear and, thus, a
     homomorphism of $A-\coring{D}$--bicomodules. On the other
     hand, if $f:M\rightarrow N$ is  homomorphism of
$A-\coring{D}$-bicomodules, then we have the commutativity of the
following diagram
\[
   \xymatrix{M\ar[r]^-{\cong}\ar[d]_-{f} &
  A\tensfun{A}M\ar[r]^-{\theta\tensfun{A}M}
   & A\tensfun{}A\tensfun{A}M\ar[r]^-{\cong}
  & A\tensfun{}M\ar[d]^-{A\tensfun{}f}
  \\ N\ar[r]_-{\cong} &
  A\tensfun{A}N\ar[r]_-{\theta\tensfun{A}N} &
  A\tensfun{}A\tensfun{A}N\ar[r]_-{\cong}
   & A\tensfun{}N.}
   \]
   Thus we have defined a natural transformation
   $\nu:1_{\bcomod{A}{\coring{D}}}\rightarrow G'F'$. Moreover,
   for every $M\in{\bcomod{A}{\coring{D}}}$, one
   has $c_{M}\nu_{M}=M$. This proves that $F'$ is separable, and
   the theorem holds.
\end{proof}

An $A$--coring $\coring{C}$ is called \emph{semisimple} if
$\coring{C}$ satisfies one of the equivalent conditions in the
following Theorem \ref{semi1}. Here, a comodule is said to be
simple if it has just two subcomodules. A semisimple comodule is a
direct sum of simple comodules.

\begin{theorem}[\cite{ElKaoutit/Gomez/Lobillo:2001unp}]\label{semi1}
Let $\coring{C}$ be an $A$--coring. The following statements are
equivalent:
\begin{enumerate}[(i)]
\item $\coring{C}$ is semisimple as a left $\coring{C}$--comodule and $\coring{C}_A$ is flat;
\item every left $\coring{C}$--comodule is semisimple and $\coring{C}_A$ is flat;
\item $\coring{C}$ is semisimple as a right $\coring{C}$--comodule and ${}_A\coring{C}$ is flat;
\item every right $\coring{C}$--comodule is semisimple and ${}_A\coring{C}$ is flat;
\item every (left or right) $\coring{C}$--comodule is semisimple, and ${}_A\coring{C}$
and $\coring{C}_A$ are projectives.
\end{enumerate}
\end{theorem}

Coseparable corings over separable algebras can be characterized
in terms of semisimplicity. The first step to obtain the
characterization is the following.

\begin{theorem}\label{cosepsemisimp}
Let $\coring{C}$ be a coseparable $A$--coring such that
${}_A\coring{C}$ is flat. If $A$ is a separable $k$--algebra, then
the $A \tensor{k} B$--coring $\coring{C} \tensor{k} \coring{D}$ is
semisimple for every semisimple $B$--coring $\coring{D}$.
\end{theorem}
\begin{proof}
Since ${}_A\coring{C}$ and ${}_B\coring{D}$ are flat, we get that
$\coring{C} \tensor{k} \coring{D}$ is flat as a left $A \tensor{k}
B$--module. By Proposition \ref{comodbicomod}, we have an
isomorphism of categories $\lcomod{\coring{C} \tensor{k}
\coring{D}} \cong \bcomod{\coring{C}}{\coring{D}^{\circ}}$. By
Theorem \ref{separaD}, the forgetful functor
$\bcomod{\coring{C}}{\coring{D}^{\circ}} \rightarrow
\rcomod{\coring{D}^{\circ}} \cong \lcomod{\coring{D}}$ is
separable. Therefore we obtain a separable functor
$\lcomod{\coring{C} \tensor{k} \coring{D}} \rightarrow
\lcomod{\coring{D}}$. Now, $\coring{D}$ is assumed to be
semisimple, which implies, by \cite[Proposition
1.2]{Nastasescu/VandenBergh/VanOystaeyen:1989}, that $\coring{C}
\tensor{k} \coring{D}$ is semisimple.
\end{proof}

\begin{remark}
The $A \tensor{k} B$--coring $\coring{C} \tensor{k} \coring{D}$
fails to be semisimple if $A$ is not assumed to be separable. For
example, the trivial structure of $A$--coring given on $A$ by the
canonical isomorphism $A \cong A \tensor{A} A$ is obviously
coseparable. Assume $k$ to be a field. If $A$ is not semisimple as
a ring, then $A \cong A \tensor{k} k$ is not semisimple as a $A
\tensor{k} k \cong A$--coring (but $k$ is obviously a semisimple
$k$--coring).
\end{remark}

Our next goal is to obtain a characterization of coseparable
corings over separable algebras over a field that generalizes the
given for separable coalgebras in \cite[Theorem
3.4]{Castano/Gomez/Nastasescu:1997}. To achieve this, we need to
prove some results which are also of independent interest.

Let us start with a generalization of the so called fundamental
theorem of coalgebras. The notion which replaces the finite
dimensionality over fields is given by the following definition.

\begin{definition}
An $R-S$--bimodule $M$ is said to be of \emph{finite type} if it
is a finitely generated left $R \tensor{k} S^{\circ}$--module,
that is, if $M = \sum_i Rm_iS$, for some $m_i \in M$, where $i$
runs over a finite set.
\end{definition}

Following \cite{Sweedler:1975}, let us endow the left dual
${}^*\coring{C} = \hom{A}{\coring{C}}{{}_AA}$ with a structure of
ring by the defining its multiplication by $gf = g(\coring{C}
\tensor{A} f)\Delta$. The map $A^{\circ} \rightarrow
{}^*\coring{C}$ induced by the counit $\epsilon: \coring{C}
\rightarrow A$ is then a ring homomorphism. Analogously, the right
dual $\coring{C}^* = \hom{A}{\coring{C}}{A_A}$ is a ring with
multiplication $gf=f(g\tensor{A} \coring{C})\Delta$, and we have a
ring homomorphism $A^{\circ} \rightarrow \coring{C}^*$. Then the
obvious structure of $A^{\circ}$--bimodule of $\coring{C}$ extends
to a ${}^*\coring{C}-\coring{C}^*$--bimodule structure. Here, the
left ${}^*\coring{C}$--action is given by $f . c = (f \tensor{A}
\coring{C})\Delta$, and the right $\coring{C}^*$--action is given
by $c.f = (\coring{C} \tensor{A} f)\Delta$.

\begin{proposition}\label{localfinie}
Let $\coring{C}$ be an $A$--coring such that ${}_A\coring{C}$ and
$\coring{C}_A$ are projective modules. Then $\coring{C} =
\cup_{i}\coring{C_i}$, where the $\coring{C}_i$'s are
$\coring{C}$--subbicomodules of $\coring{C}$ which are of finite
type as $A$--bimodules.
\end{proposition}
\begin{proof}
 Thinking of
$\coring{C}$ as a left ${}^*\coring{C} \tensor{k}
{\coring{C}^*}^{\circ}$--module, we can write it as a direct union
of its finitely generated submodules. Thus, $\coring{C} =
\cup_{i}\coring{C}_i$, where the $\coring{C}_i$ are
${}^*\coring{C}-\coring{C}^*$--subbimodules of $\coring{C}$ of
finite type. Let us prove that each $\coring{C}_i$ is of finite
type as an $A$--bimodule. By \cite[Corollary
2.11]{ElKaoutit/Gomez/Lobillo:2001unp}, the
$\coring{C}$--subbicomodules of $\coring{C}$ are precisely its
${}^*\coring{C}-\coring{C}^*$--subbimodules. Thus, accordingly
with \cite[Proposition 2.10]{ElKaoutit/Gomez/Lobillo:2001unp},
this proof is finished as soon as we prove that every rational
${}^*\coring{C}-\coring{C}^*$--bimodule of finite type $M$ is of
finite type as an $A$--bimodule. For, if $M = \sum_s
{}^*\coring{C} m_s \coring{C}^*$, for finitely many $m_s \in M$,
then we proceed as follows: choose a system of left rational
parameters $(c_{sj},m_{sj})$ for each $m_s$ (see
\cite{ElKaoutit/Gomez/Lobillo:2001unp} for this notion), and then
choose a system of right rational parameters $(m_{sjk},c_{sjk})$
for every $m_{sj}$. If $m \in M$, then $m = \sum_s b_s m_s b'_s$,
for some $b_s \in {}^*\coring{C}, b'_s \in \coring{C}^*$. Compute
\[
m=\sum_{s,j}(m_{sj}b_{s}(c_{sj}))b'_{s}=\sum_{s,j}(m_{sj}b'_{s})b_{s}(c_{sj})=
   \sum_{s,j,k}b'_{s}(c_{sjk})m_{sjk}b_{s}(c_{sj}),
\]
whence $M = \sum_{s,j,k}Am_{sjk}A$, and it is an $A$--bimodule of
finite type. This finishes our proof.
\end{proof}

A coring $\coring{C}$ is said to be \emph{simple} if it has no
nontrivial subbicomodules. Examples of simple $A$--corings are the
Sweedler's canonical corings $A \tensor{B} A$, where $B$ is a
simple artinian subring of $A$ (see \cite[Proposition 4.2
]{ElKaoutit/Gomez/Lobillo:2001unp}).

\begin{corollary}\label{simplefinito}
If $\coring{C}$ is a simple coring such that ${}_A\coring{C}$ and
$\coring{C}_A$ are projective modules, then $\coring{C}$ is an
$A$--bimodule of finite type.
\end{corollary}

Every homomorphism of $A$--corings $\varphi : \coring{C}
\rightarrow \coring{D}$ induces a functor $(-)_{\varphi} :
\lcomod{\coring{C}} \rightarrow \lcomod{\coring{D}}$, called
corestriction functor. The corresponding homomorphism of rings
$\varphi^* : \coring{D}^* \rightarrow \coring{C}^*$ (see
\cite[Proposition 3.2]{Sweedler:1975}) gives the usual restriction
of scalars functor $(-)_{\varphi^*} : \rmod{\coring{C}^*}
\rightarrow \rmod{\coring{D}^*}$.

\begin{proposition}\label{cosepsep}
Let $\varphi : \coring{C} \rightarrow \coring{D}$ be a
homomorphism of $A$--corings. If $\coring{C}_A$ and
$\coring{D}_A$ are finitely generated projective modules, then
the following conditions are equivalent:
\begin{enumerate}[(i)]
\item The functor $(-)_{\varphi} : \lcomod{\coring{C}} \rightarrow
\lcomod{\coring{D}}$ is separable;
\item the functor $(-)_{\varphi^*} : \rmod{\coring{C}^*}
\rightarrow \rmod{\coring{D}^*}$ is separable.
\end{enumerate}
\end{proposition}
\begin{proof}
By \cite[Lemma 4.3]{Brzezinski:2000unp}, we have isomorphisms of
categories $\lcomod{\coring{C}}\cong\rmod{\coring{C}^*}$ and
    $\lcomod{\coring{D}}\cong\rmod{\coring{D}^*}$. We get the result from the
    commutativity of the following diagram of functors
    \[\xymatrix{\lcomod{\coring{C}}\ar[r]^-{(-)_\varphi}\ar[d]^-{\cong} &
     \lcomod{\coring{D}}\ar[d]^-{\cong} \\
     \rmod{\coring{C}^*}\ar[r]_-{(-)_{\varphi^*}} &
     \rmod{\coring{D}^*}.}
     \]
\end{proof}

\begin{theorem}\label{cosepsepalgebra}
Let $\coring{C}$ be a coring over a separable $k$--algebra $A$.
If $\coring{C}_A$ is projective and finitely generated then
$\coring{C}$ is a coseparable coring if and only if
$\coring{C}^*$ is a separable $k$--algebra.
\end{theorem}
\begin{proof}
Consider the counit $\epsilon : \coring{C} \rightarrow A$ as a
homomorphism of corings. Then the functor $(-)_{\epsilon} :
\lcomod{\coring{C}} \rightarrow \lmod{A}$ is just the forgetful
functor, and, in this way, $(-)_{\epsilon}$ is separable if and
only if $\coring{C}$ is a coseparable coring. On the other hand,
the functor $(-)_{\epsilon^*}: \rmod{\coring{C}^*} \rightarrow
\rmod{A^{\circ}}$ is nothing but the restriction of scalars
functor associated to the ring homomorphism $\epsilon^* :
A^{\circ} \rightarrow \coring{C}^*$. Thus, by Proposition
\ref{cosepsep}, $\coring{C}$ is coseparable if and only if
$\rmod{\coring{C}^*} \rightarrow \rmod{A^{\circ}}$ is separable.
On the other hand, by \cite[Proposition 1.3
]{Nastasescu/VandenBergh/VanOystaeyen:1989}, we have that the
restriction of scalars functor $\rmod{A^{\circ}} \rightarrow
\rmod{k}$ is separable, since $A^{\circ}$ is a separable
$k$--algebra. Now, the composite $\rmod{\coring{C}^*} \rightarrow
\rmod{A^{\circ}} \rightarrow \rmod{k}$ is the evident restriction
of scalars functor. Hence, by \cite[Lemma
1.1]{Nastasescu/VandenBergh/VanOystaeyen:1989}, this functor is
separable if and only if $\rmod{\coring{C}^*} \rightarrow
\rmod{A^{\circ}}$ is. In other words, $\coring{C}$ is a
coseparable coring if and only if $\coring{C}^*$ is a separable
$k$--algebra.
\end{proof}

Let $R$ be a ring and $M$ a right $R$--module. Recall that an
$R$--submodule $N$ of $M$ is said to be \emph{pure} if the
canonical map $N \tensor{R} L \rightarrow M \tensor{R} L$ is
injective for every left $R$--module $L$. Pure submodules of left
$R$--modules are similarly defined. We say that a homomorphism of
right $R$--modules $f : X \rightarrow Y$ is pure if $Im f$ is a
pure $R$--submodule of $Y$.

\begin{proposition}\label{subsep}
Let $\varphi : \coring{C} \rightarrow \coring{D}$ be a
homomorphism of $A$--corings. If $\varphi$ is a pure monomorphism
of right $A$--modules, then the functor $(-)_{\varphi} :
\lcomod{\coring{C}} \rightarrow \lcomod{\coring{D}}$ is separable.
\end{proposition}
\begin{proof}
Let $M,N\in\lcomod{\coring{C}}$. If $f:M\rightarrow N$
    is a homomorphism of $\coring{D}$-comodules, then it is
    a homomorphism of $\coring{C}$-comodules. In fact, one has
    \[
    (\varphi\tensor{A}N)\lambda_{N}f=
    (\coring{D}\tensor{A}f)(\varphi\tensor{A}M)\lambda_{M}=
    (\varphi\tensor{A}N)(\coring{C}\tensor{A}f)\lambda_{M}.
    \]
    Since $\varphi$ is a pure monomorphism of right $A$--modules,
    we get that
   $\lambda_{N}f=(\coring{C}\tensor{A}f)\lambda_{M}$. Therefore we
   can define the map
   $\varphi_{M,N}:\hom{\coring{D}}{M}{N}\rightarrow\hom{\coring{C}}{M}{N}$
   by $\varphi_{M,N}(f)=f$. Thus, by definition
   (see \cite[p. 398]{Nastasescu/VandenBergh/VanOystaeyen:1989}),
   $(-)_{\varphi}$ is a separable functor.
\end{proof}

A subbicomodule $\coring{D}$ of a coring $\coring{C}$ is said to
be a \emph{subcoring} if $\coring{D}$ is pure in $\coring{C}$ both
as a left and as a right $A$--submodule. In such a case, the
restriction of the comultiplication of $\coring{C}$ to
$\coring{D}$ gives a well defined structure of $A$--coring on
$\coring{D}$ such that the inclusion $\coring{D} \subseteq
\coring{C}$ is a homomorphism of $A$--corings. Every subbicomodule
of a semisimple coring is a subcoring, since it is a direct
summand as a left and as a right $A$--submodule. The following is
a consequence of Proposition \ref{subsep} and \cite[Proposition
1.2]{Nastasescu/VandenBergh/VanOystaeyen:1989}.

\begin{corollary}\label{semisimpfinsemisimp}
Every subcoring of a semisimple coring is semisimple.
\end{corollary}

The following result tell us that the separability of a morphism
of corings defined on a semisimple coring is a `local' property.

\begin{proposition}\label{seplocal}
Let $\coring{C}$ be a semisimple $A$--coring and let $\varphi :
\coring{C} \rightarrow \coring{D}$ be a homomorphism of
$A$--corings such that the image of $\varphi$ is pure in
$\coring{D}$ both as a left and right $A$--module. Then
$(-)_{\varphi} : \lcomod{\coring{C}} \rightarrow
\lcomod{\coring{D}}$ is a separable functor if and only if the
functor $(-)_{\varphi|\coring{C}'} : \lcomod{\coring{C}'}
\rightarrow \lcomod{\coring{D}'}$ is separable for every pair of
subcorings $\coring{C}' \leq \coring{C}$ and $\coring{D'} \leq
\coring{D}$ of finite type as $A$--bimodules such that
$\varphi(\coring{C}') \subseteq \coring{D}'$.
\end{proposition}
\begin{proof}
Assume $(-)_{\varphi}$ to be a separable functor, and let
$\coring{C}'$ and $\coring{D}'$ as declared. If $i : \coring{C}'
\rightarrow \coring{C}$ and $j:\coring{D}' \rightarrow \coring{D}$
denote the inclusion maps, then $(-)_j \circ
(-)_{\varphi|\coring{C}'} = (-)_{\varphi} \circ (-)_i$, which is
separable by Proposition \ref{subsep}. Apply \cite[Lemma 1.1
]{Nastasescu/VandenBergh/VanOystaeyen:1989} to obtain that
$(-)_{\varphi|\coring{C}'}$ is separable. \\
Conversely,
 since $\coring{C}$ is
semisimple, we get from \cite[Theorem
3.5]{ElKaoutit/Gomez/Lobillo:2001unp} a decomposition $\coring{C}
= \oplus_i \coring{C}_i$, with $\coring{C}_i$ a simple semisimple
$A$--coring for every $i$. By \ref{simplefinito}, all the
$\coring{C}_i$'s are $A$--bimodules of finite type. Write
$\coring{D}_i = \varphi(\coring{C}_i)$, which is pure in
$\coring{D}$ as a left and as a right $A$--submodule for every
$i$. Therefore, these $\coring{D}_i$'s are subcorings of
$\coring{D}$ which are of finite type as $A$--bimodules. Let us
denote by $\varphi_i : \coring{C}_i \rightarrow \coring{D}_i$ the
coring homomorphisms induced by $\varphi$. By hypothesis, the
functors $(-)_{\varphi_i}$ are separable which implies, by
\cite[Theorem 4.7]{Gomez:2001unp}, that there exists bicomodule
homomorphisms $\omega_i : \coring{C}_i \cotensor{\coring{D}_i}
\coring{C}_i \rightarrow \coring{C}_i$ such that
$\omega_{\coring{C}_i}\overline{\Delta_{\coring{C}_i}} =
\coring{C}_i$ for every $i$. Here,
$\overline{\Delta}_{\coring{C_i}}$ is given by the
comultiplication on $\coring{C}_i$. Finally,
$\coring{C}_i\cotensor{\coring{D}_i}\coring{C}_i=
     \coring{C}_i\cotensor{\coring{D}}\coring{C}_i$ and
     $\coring{C}\cotensor{\coring{D}}\coring{C}=\bigoplus_i
     (\coring{C}_i\cotensor{\coring{D}}\coring{C}_i)\oplus
     \bigoplus_{i\neq
     j}\coring{C}_i\cotensor{\coring{D}}\coring{C}_j$.
      Define $\omega=\bigoplus_i\omega_i$ over
      $\bigoplus_i(\coring{C}_i\cotensor{\coring{D}}\coring{C}_i)$
      and $ \omega=0$ over
      $\bigoplus_{i\neq
      j}\coring{C}_i\cotensor{\coring{D}}\coring{C}_j$;
      we have $\omega\overline{\Delta}=\coring{C}$. By
      \cite[Theorem 4.7]{Gomez:2001unp}, $(-)_{\varphi}$ is a
      separable functor.
\end{proof}

If we apply Proposition \ref{seplocal} to the homomorphism of
corings given by the counity $\epsilon : \coring{C} \rightarrow
A$, then we obtain the following.

\begin{corollary}\label{cosepfincosep}
A semisimple $A$--coring is coseparable if and only every
subcoring of finite type as $A$--bimodule is coseparable.
\end{corollary}

Let $\coring{C}$ be an $A$--coring and $K$ a commutative ring
containing $k$. Then $\coring{C} \tensor{k} K$ is a coring over
$A \tensor{k} K$, with left dual ring ${}^*(\coring{C} \tensor{k}
K) =  \hom{A \tensor{k} K}{\coring{C} \tensor{k} K}{A \tensor{k}
K}$ endowed with the convolution product. On the other hand, we
have the tensor product structure of ring on ${}^*\coring{C}
\tensor{k} K$. These two $k$--algebras are related.

\begin{proposition}\label{dualKdual}
The homomorphism of $k$-modules
\[\Psi:\hom{A}{\coring{C}}{A}\tensor{k}K
   \rightarrow\hom{A\tensor{k}K}{\coring{C}\tensor{k}K}{A\tensor{k}K}
   \] defined by
   $\Psi(f\tensor{k}\alpha)(x\tensor{k}\beta)=f(x)\tensor{k}\alpha\beta$
    for every $f\in\hom{A}{\coring{C}}{A}\tensor{k}K,\ x\in\coring{C},\ \alpha,\ \beta\in K$
    is a homomorphism of $k$--algebras. Moreover, if
    ${}_A\coring{C}$ is a finitely generated projective module, then $\Psi$ is an isomorphism.
\end{proposition}
\begin{proof}
Let us prove that $\Psi$ is a homomorphism of algebras. Let
$f,f'\in\hom{A}{\coring{C}}{A}, x\in\coring{C}$ and
$\alpha,\alpha',\beta\in K$. If $\Delta(x)=x_1\tensor{k}x_2$, then
\[
\begin{split}
\Psi[(f\tensor{k}\alpha)(f'\tensor{k}\alpha')](x\tensor{k}\beta)
&= ff'(x)\tensor{k}\alpha\alpha'\beta \\ &=
 f(x_1f'(x_2))\tensor{k}\alpha\alpha'\beta \\ &=
\Psi(f\tensor{k}\alpha)[x_1f'(x_2)\tensor{k}\beta\alpha']
\\ &=
 \Psi(f\tensor{k}\alpha)[(x_1\tensor{k}\beta)\Psi(f'\tensor{k}\alpha')(x_2\tensor{k}1_k)]
 \\ &=
 [\Psi(f\tensor{k}\alpha)\Psi(f'\tensor{k}\alpha')](x\tensor{k}\beta).
 \end{split}
 \]
Now assume that ${}_A\coring{C}$ is a finitely generated
projective module. Let
$\Phi:\hom{A\tensor{k}K}{\coring{C}\tensor{k}K}{A\tensor{k}K}
    \rightarrow\hom{A}{\coring{C}}{A}\tensor{k}K$
   the map defined as follows: First, take a dual basis $(x_i,\varphi_i)_i$ for ${}_A\coring{C}$.
   If
   $g\in\hom{A\tensor{k}K}{\coring{C}\tensor{k}K}{A\tensor{k}K}$,
   write
   $g(x_i\tensor{k}1)=\sum_i(a_{ij}\tensor{k}\alpha_{ij})$, and
   define $\Phi(g)=\sum_{i,j}\varphi_ia_{ij}\tensor{k}\alpha_{ij}$. It is
   not hard to see that $\Phi$ is the inverse map for $\Psi$.
\end{proof}

We are now in position to prove our main result.

\begin{theorem}
Let $A$ be a separable algebra over a field $k$. The following
conditions are equivalent for an $A$--coring $\coring{C}$.
\begin{enumerate}[(i)]
\item $\coring{C}$ is coseparable;
\item $\coring{C} \tensor{k} K$ is a semisimple $A \tensor{k} K$--coring for every field
extension $k \subseteq K$;
\item $\coring{C} \tensor{k} \coring{D}$ is a semisimple $A
\tensor{k} B$--coring for every semisimple $B$--coring
$\coring{D}$ ($B$ is any $k$--algebra);
\item $\coring{C} \tensor{k} \coring{C}^{\circ}$ is a semisimple $A \tensor{k} A^{\circ}$--coring.
\end{enumerate}
\end{theorem}
\begin{proof}
First, we recall that our separable algebra $A$ over the field $k$
is necessarily a finite dimensional semisimple $k$--algebra.\\
$(i) \Rightarrow (iii)$ This is a particular case of Theorem
\ref{cosepsemisimp}.\\
$(iii) \Rightarrow (ii)$ Put $B = \coring{D} = K$.\\
$(ii) \Rightarrow (i)$ The coring $\coring{C}$ is semisimple, so
by Corollary \ref{cosepfincosep} and Corollary
\ref{semisimpfinsemisimp} we can assume that $\coring{C}$ is an
$A$--bimodule of finite type. This implies that $\coring{C}$ is
finite dimensional as a $k$--vector space and, in particular,
${}_A\coring{C}$ is finitely generated and projective. By
Proposition \ref{dualKdual}, ${}^*\coring{C} \tensor{k} K \cong
{}^*(\coring{C} \tensor{k} K)$ for every field extension $k
\subseteq K$. This last $k$--algebra is semisimple because the $A
\tensor{k} K$--coring $\coring{C} \tensor{k} K$ is semisimple, and
$\coring{C} \tensor{k} K$ is finitely generated and projective as
a left $A \tensor{k} K$--module. By a classical result
(\cite[Theorem 2.5]{DeMeyer/Ingraham:1971}), ${}^*\coring{C}$
becomes a separable $k$--algebra which implies, by Theorem
\ref{cosepsepalgebra},
that $\coring{C}$ is a coseparable coring.\\
$(iii) \Rightarrow (iv)$ Taking $\coring{D} = B = k$ we obtain
that $\coring{C}$ is semisimple. By Theorem \ref{semi1},
$\coring{C}^{\circ}$ is
also semisimple. Therefore, $\coring{C} \tensor{k} \coring{C}^{\circ}$ is semisimple.\\
$(iv) \Rightarrow (i)$ The comultiplication $\Delta : \coring{C} \rightarrow \coring{C}
 \tensor{A} \coring{C}$ is a homomorphism of $\coring{C}$--bicomodules.
 By Proposition \ref{comodbicomod}, it is a homomorphism of left
 $\coring{C} \tensor{k} \coring{C}^{\circ}$--comodules. Since this
 coring is assumed to be semisimple, $\Delta$ splits as a left
 $\coring{C} \tensor{k} \coring{C}^{\circ}$--comodule map and, hence, as a
 $\coring{C}$--bicomodule map. Therefore, $\coring{C}$ is
 coseparable.
\end{proof}

\providecommand{\bysame}{\leavevmode\hbox
to3em{\hrulefill}\thinspace}
\providecommand{\MR}{\relax\ifhmode\unskip\space\fi MR }
\providecommand{\MRhref}[2]{%
  \href{http://www.ams.org/mathscinet-getitem?mr=#1}{#2}
} \providecommand{\href}[2]{#2}

\end{document}